\documentclass[12pt]{article}
\usepackage{amsbsy}
\usepackage{subfigure,amsmath,amssymb}
\usepackage{float}
\usepackage{mathrsfs}
\usepackage{amssymb,url}
\usepackage{cite}
\usepackage{graphicx}
\usepackage[]{caption2}

\usepackage[colorlinks,linkcolor=red,urlcolor=blue,linktocpage]{hyperref}
\numberwithin{equation}{section}

\begin{document}
\title{\Large \bf Uniqueness of solutions to the logarithmic Minkowski problem in $\mathbb{R}^3$
}
\author{ \small \bf Shibing Chen~~\bf Yibin Feng~~\bf Weiru Liu
\\ \small  School of Mathematical Sciences, University of Science and Technology of China,\\
\small Hefei, 230026, China \\  \small E-mails:
 chenshib@ustc.edu.cn, fybt1894@ustc.edu.cn, \\  \small lwr19997@mail.ustc.edu.cn
}

\date{}
\maketitle

\vskip 20pt

\begin{center}
\begin{minipage}{12cm}
\small
 {\bf Abstract:} 
 %Andrews (Invent. Math. 138, 1999) proved the uniqueness of solutions to the logarithmic Minkowski problem in $\mathbb{R}^3$ when the measure $\mu$ has a constant density.
  % In $\mathbb{R}^3$, this paper presents the uniqueness of solutions to the problem when $\mu$ has the density (not necessarily constant) close to $1$.
 In this paper, we prove the uniqueness of solutions to the logarithmic Minkowski problem in $\mathbb{R}^3$ without symmetry condition, provided the density of the measure is close to $1$ in $C^{\alpha}$ norm. This result also implies the uniqueness of self-similar solutions to the anisotropic Gauss curvature flow in $\mathbb{R}^3$ when the speed function is $C^{\alpha}$ close to a positive constant.

% {\bf Keywords:} Cone-volume measure, Minkowski problem, logarithmic Minkowski problem, Monge-Amp\`{e}re equation. \\
% {\bf 2010 Mathematics Subject Classification:} 52A40.

 \vskip 0.1cm
\end{minipage}
\end{center}

\vskip 20pt

\section{\bf Introduction}

 \ \ \ \ Let $\mathbb{R}^n$ be the $n$-dimensional Euclidean space. The unit sphere in $\mathbb{R}^n$ is denoted by $\mathbb{S}^{n-1}$. A convex body in $\mathbb{R}^n$ is a compact convex set with non-empty interior. Denote by $\mathcal{K}^n$ the set of all convex bodies in $\mathbb{R}^n$ that contain the origin in their closure. For a convex body $K$, its support function $h_K: \mathbb{S}^{n-1}\rightarrow \mathbb{R}$ is defined by 
 $h_K(x)=\max\{x\cdot y: y\in K\}$, where ``$\cdot$"  is the standard inner product in $\mathbb{R}^n$.
 
 The cone-volume measure $V_K$ of a convex body $K\in \mathcal{K}^n$ is a Borel measure on $\mathbb{S}^{n-1}$, defined for a Borel $\omega\subset \mathbb{S}^{n-1}$ by
 \begin{eqnarray}\label{1.1}
 	V_K(\omega)=\frac{1}{n}\int_{x\in \nu_K^{-1}(\omega)}x\cdot \nu_K(x)d\mathcal{H}^{n-1}(x),
 \end{eqnarray}
where $\nu_K:\partial K\rightarrow \mathbb{S}^{n-1}$ is the set valued normal mapping, namely,
$\nu_K(x):=\{\nu\in\mathbb{S}^{n-1}: \nu\cdot x\geq \nu\cdot z\ \text{for any}\ z\in K\},$  
% the set of boundary points of $K$ that have a unique outer unit normal,
 and $\mathcal{H}^{n-1}$ is $(n-1)$-dimensional Hausdorff measure.  Indeed, $V_K(\omega)$ is the volume of the cone constituted of the segments connecting the origin and $\nu_K^{-1}(\omega).$
The cone-volume measure has clear geometric significance and has attracted great attention from many scholars; see for example \cite{12a, 14a, 15a, 16a, 25a, 26a, 33a, HLL, Xiong}.

One of the cornerstones of convex geometry analysis is the Minkowski problem. It is the characterisation problem of surface area measure of convex bodies, and was first studied by Minkowski himself. After Minkowski's original work, the associated Minkowski type problems have been extensively investigated (see e.g. \cite{27a, 37ab, 39a, 46a, 26a, 15a, 44a, 45a, 40a, 25a, 41a, 42a, 1b, 2b, 3b, 4b, 5b, 6b, 7b, 9b, 10b, 11b, 12b, 13b, 14b, 15b, 29a, 43a, WXL, GL}), where the logarithmic Minkowski problem is one of the most central Minkowski type problems and is the problem of characterizing the cone-volume measure.

\noindent{\bf Logarithmic Minkowski problem.} Find necessary and sufficient conditions on a finite Borel measure $\mu$ on $\mathbb{S}^{n-1}$ so that $\mu$ is the cone-volume measure of a convex body $K$ in $\mathbb{R}^n$.

When $n\mu=f\mathcal{H}^{n-1}$, the associated partial differential equation for the logarithmic Minkowski problem is the following Monge-Amp\`{e}re type equation on $\mathbb{S}^{n-1}$:
\begin{eqnarray}\label{1.2}
h_K\mbox{det}(\nabla^2h_K+h_KI)=f,
\end{eqnarray}
where $\nabla^2h_K=\left(\nabla_{ij}h_K\right)$ is the Hessian matrix of covariant derivatives of $h_K$ with respect to an orthonormal frame on $\mathbb{S}^{n-1}$, and $I$ is the identity matrix. The solutions of \eqref{1.2} can also be interpreted as the self-similar solutions to the anisotropic Gauss curvature flow
\begin{equation}\label{gcf}
u_t(x, t)=f(\nu(x, t))\mathcal{K}(x, t)\nu(x, t),
\end{equation}
where $u(\cdot, t):\mathbb{S}^{n-1}\rightarrow \mathbb{R}$ is a family of time dependent embedding, namely, the evolving hypersurface,  and  $\mathcal{K}(x, t)$ (resp. $\nu(x, t)$) is the Gauss curvature (resp. unit inner normal) of $u(\mathbb{S}^{n-1}, t)$ at $u(x, t).$

The existence of solutions to the logarithmic Minkowski problem for symmetric measure has been completely solved by B\"{o}r\"{o}czky, Lutwak, Yang and Zhang \cite{15a}. The necessary and sufficient condition found in \cite{15a} also turns out to be a sufficient condition for the existence of solutions when the measure is not assumed to be symmetric.
A sufficient condition of the existence of solutions to the problem for discrete measure (not necessarily symmetric) was given by Zhu \cite{25a}, and later by B\"{o}r\"{o}czky, Heged\H{u}s and Zhu \cite{26a}.
 For smooth measure, we refer to  \cite{58a} and references therein.

 The uniqueness of solutions to the logarithmic Minkowski problem appears to be more complicated than its existence. 
 Indeed, even in dimension 2, the uniqueness of solutions to equation \eqref{1.2} fails for some positive smooth $f$, and see the example in \cite{Y06}.
  The uniqueness of solutions to the planar logarithmic Minkowski problem was studied by Gage \cite{60a} and Stancu \cite{58a, 59a} for smooth and discrete measures,  by B\"{o}r\"{o}czky, Lutwak, Yang and Zhang \cite{16a} for symmetric measures, and by  Xi and Leng \cite{X} for general measures when the convex bodies are at a dilation position.
  When $n=3$ and $\mu$ has a constant density, the uniqueness was shown by Firey \cite{61a} for symmetric measures and by Andrews \cite{48a} for non-symmetric measures. 
  In 2017,  Brendle, Choi and Daskalopoulos \cite{BCD} extended Andrews' result to general dimensions. 
  For general positive smooth $f,$ the uniqueness of solutions to \eqref{1.2} in higher dimensions is a very challenging problem.  When $f$ is an even function and close to 1 in $C^{\alpha}$ norm, Chen, Huang, Li and Liu \cite{14b}
  proved the uniqueness of solutions to \eqref{1.2} in the class of origin-symmetric convex bodies, based on the  local results in \cite{KM17, CLM}.

 In this paper, we establish the following uniqueness result when $n=3$ and $\mu$ has a positive density of $C^\alpha$ close to 1. Note that we do not require any symmetry condition in the following theorem.
 
{\it \noindent{\bf Theorem 1.1.}~Suppose $1/\lambda<f<\lambda$ for some $\lambda>0.$ Then there exists a small constant $\epsilon_0>0$ depending only on $\lambda,$ such that if $\|f-1\|_{C^\alpha}\leq \epsilon_0$ then equation (\ref{1.2}) has a unique solution in $\mathbb{R}^3$.}

{\it \noindent{\bf Remark 1.} The above result implies the uniqueness of self-similar solutions to the anisotropic Gauss curvature flow \eqref{gcf} in dimension 3. 
}

The main ingredient in the proof of the above theorem is a uniform upper bound for  solutions to equation \eqref{1.2}.
%see the statement of Lemma 3.1. 

{\it \noindent{\bf Lemma 1.}
Suppose $1/\lambda<f<\lambda$ for some $\lambda>0$ and $K\in \mathcal{K}^n$ in $\mathbb{R}^3$. If $K$ is a solution to $h_KdS_K=fd\mathcal{H}^{2}$, then $\|h_K\|_{L^\infty(\mathbb{S}^{2})}\leq C$ for some constant $C$ depending only on $\lambda$.}

When the dimension $n=2,$ the estimate similar to Lemma 1 is available in the literature \cite{13b}. Indeed, a
stronger result was proved there, namely, the authors also established a uniform positive lower bound of solutions to \eqref{1.2}.
However, in dimension higher than 2, things are quite different. Actually,
there exist examples showing that 
\eqref{1.2} has no uniform positive lower bound, namely, there exists positive smooth $f$ such that the resulting convex body may have the origin on its boundary, see \cite{46a} for more details. The upper bound is also a challenging question for \eqref{1.2} in dimension $n\geq 3.$  Note that the proof in dimension 2 highly relies on the ODE structure of \eqref{1.2}, and the maximum principle turns out to be not very useful to establish the estimate.

In this paper, we develop a novel blow-down argument to prove the uniform upper bound.  In dimension 3, if there is no uniform upper bound for  solutions of \eqref{1.2}, then by  the geometric meaning of the equation we can find a sequence of solutions  of \eqref{1.2} with righthand side $f$ bounded from below and above by some fixed positive constant, such that the convex bodies corresponding to this sequence  either converge to an infinite straight line, or converge to a two dimensional hyperplane.

The former case can be ruled out by choosing proper testing domains on $\mathbb{S}^{n-1}$ and using the definition of cone volume measure. The latter case is more delicate. By blowing down, up to a subsequence it converges to a bounded two dimensional convex set containing the origin on its boundary.
Then a contradiction is made by using very delicate convex analysis technique. Roughly speaking, we will show that the cone volume measure on $\mathbb{S}^1$, induced by
the two dimensional blow-down set,
%inherits from the original sequence of 3 dimensional convex \textcolor{red}{ bodies and the }property that t 
has density bounded from below and above, which forces the origin to be away from the boundary. This contradicts the property that the blow-down set contains the origin on its boundary. 
Once the uniform bound is established, Theorem 1.1 follows from the regularity theory of Monge-Amp\`ere equation and the inverse function theorem.

{\it \noindent{\bf Remark 2.}
If the dimension $n>3,$  then the blow-down set may be 3 dimensional. However, when $n=3,$ a solution to \eqref{1.2} may attain 0, namely the corresponding convex body may contain the origin on its boundary, even assuming $f$ to be positive smooth. Hence we can not reach a contradiction as above. It would be interesting to see whether Lemma 1 holds or not in higher dimensions.  
}
%Let $K$ be the convex body corresponding to the solution of \eqref{1.2} in dimension 3.
 %By John's Lemma we have that $\frac{1}{8}E\subset K\subset 8E$ for some ellipsoid $E\subset \mathbb{R}^3$ centered at the center of mass of $E$. By a rotation of coordinates we may assume the principal radii of $E$ are $r_ie_i,\ i=1,2,3,$ where $e_i$ are coordinate directions and $r_1\leq r_2\leq r_3.$ Observe that 
 %$1/C\leq vol(E)\approx\ vol(K)\leq C$ for some positive constant $C$ depending only on the lower and upper bound of $f.$ If the solutions of \eqref{1.2} have no uniform upper bound, then there exists solution $K$ such that $r_1<<1$ and $r_3>>1.$ There are two cases to be ruled out: case 1, $r_1

The rest of the paper is organized as follows. In section 2, we list some notions and basic facts regarding convex bodies.  In section 3, we prove the uniform upper bound of solutions to \eqref{1.2}.
The proof of Theorem 1.1 is presented in section 4.

 \section{\bf Preliminaries}

 \ \ \ \ \ In this section, we recall some notions and results in the theory of convex bodies, and for more details we refer the reader to Gardner \cite{62a} and Schneider \cite{63a}.

We write $|x|=\sqrt{x\cdot x}$ for $x\in \mathbb{R}^n$ and $B_r(x)=\{y\in \mathbb{R}^n: |x-y|\leq r\}$ for the Euclidean ball in $\mathbb{R}^n$ with center $x\in \mathbb{R}^n$ and radius $r>0$. Let $o$ be the origin in $\mathbb{R}^n$. When $x=o$ and $r=1$, $B_1(o)$ is usually abbreviated as $B$.

The Hausdorff distance of two convex bodies $K, L$ in $\mathbb{R}^n$ is defined by
\begin{eqnarray*}
d(K, L)=\max_{u\in \mathbb{S}^{n-1}}\left|h_K(u)-h_L(u)\right|.
\end{eqnarray*}
Let $K_i$ be a sequence of convex bodies. We say $K_i$ converges to a convex body $K$ if
$d(K_i, K)\rightarrow 0$ as $i\rightarrow \infty$.

The surface area measure, $S_K$, of a convex body $K$ is a Borel measure on $\mathbb{S}^{n-1}$, which is defined for a Borel $\omega\subset \mathbb{S}^{n-1}$ by
\begin{eqnarray}\label{2.1}
S_K(\omega)=\mathcal{H}^{n-1}(\nu_K^{-1}(\omega)).
\end{eqnarray}
For convenience, 
in the following we will regard $\nu_K:\partial K\rightarrow \mathbb{S}^{n-1}$ as a set valued function, namely,
$\nu_K(x):=\{\nu\in\mathbb{S}^{n-1}: \nu\cdot x\geq \nu\cdot z\ \text{for any}\ z\in K\}.$  

We will make use of the weak continuity of surface area measure, i.e., if $K_i$ is a sequence of convex bodies and $K$ is a convex body then
\begin{eqnarray*}
\lim_{i\rightarrow \infty}K_i=K \Longrightarrow \lim_{i\rightarrow \infty}S_{K_i}=S_K, \mbox{weakly}.
\end{eqnarray*}
Let $V(K)$ denote the volume of a convex body $K\in \mathcal{K}^n$. Then by (\ref{1.1}) and (\ref{2.1}), we easily see
\begin{eqnarray*}
	V_K(\mathbb{S}^{n-1})=V(K)=\frac{1}{n}\int_{u\in \mathbb{S}^{n-1}}h_K(u)dS_K(u).
\end{eqnarray*}

For convex bodies $K, L$ and $x\in\mathbb{R}^n$, define
\begin{eqnarray*}
\mbox{dist}(x, K)=\inf_{y\in K}|x-y|~~\mbox{and}~~\mbox{dist}(K, L)=\inf_{x\in K, y\in L}|x-y|.
\end{eqnarray*}

 In the following sections we use the notation $a\gtrsim b$ (resp. $a\lesssim b$) if there exists a constant $C$ 
 depending only on the upper and lower bounds of $f$
 such that $a \geq C b$ (resp. $a\leq C b$), and $a\approx b$ means that $C^{-1} b\le a\le C b$.

\section{\bf A priori estimate}

\ \ \ \ \  In this section, we prove Lemma 1, namely,  we will establish a uniform upper bound of solutions to the logarithmic Minkowski problem, which is crucial for the proof of Theorem 1.1.

%{\it \noindent{\bf Lemma 3.1.}~
%Suppose $1/\lambda<f<\lambda$ for some $\lambda>0$ and $K\in \mathcal{K}^n$ in $\mathbb{R}^3$. If $K$ is a solution to $h_KdS_K=fd\mathcal{H}^{2}$, then $\|h_K\|_{L^\infty(\mathbb{S}^{2})}\leq C$ for some constant $C$ depending only on $\lambda$.}

{\bf Proof of Lemma 1.} Suppose that $K\in \mathcal{K}^n$ is a solution to the equation $h_KdS_K=fd\mathcal{H}^{2}$. By the assumption, we see
\begin{eqnarray*}
V(K)=\frac{1}{3}\int_{\mathbb{S}^2}f\approx1.
\end{eqnarray*}
%Let $E$ be the John's ellipsoid of $K,$ namely,
By John's Lemma (see \cite{d}), there exists an ellipsoid $E\subset \mathbb{R}^3$ such that
\begin{eqnarray}\label{john}
E\subset K\subset3^{\frac{3}{2}}E,
\end{eqnarray}
where $E$ is centred at the centre of mass of $K$. By a rotation of coordinates we may assume that $r_1, r_2, r_3$ are the principal radii of $E$ and the corresponding directions are $e_1, e_2, e_3$. Without loss of generality, we may assume that
\begin{eqnarray*}
r_1\leq r_2\leq r_3.
\end{eqnarray*}
We need to show
\begin{eqnarray*}
r_1\gtrsim r_3,
\end{eqnarray*}
namely
\begin{eqnarray*}
r_1\geq C r_3
\end{eqnarray*}
for some constant $C$ depending only on $\lambda$. If not, at least one of  the following two cases will happen.\\
Case I: for arbitrary large $M>0,$ there exists a convex body $K$ such that its cone volume measure has density $f_K$ satisfying  \begin{equation}\label{case1}
 1/\lambda<f_K<\lambda,\ \text{and that}\ \frac{r_3}{r_2}>M;
 \end{equation}

\noindent Case II: there exists a constant $C_2>0,$  and for any $M>1$ large, there exists a convex body $K$
such that its cone volume measure has density $f_K$ satisfying 
\begin{equation}\label{case222}
 1/\lambda<f_K<\lambda,\ \text{that}\frac{r_2}{r_1}>M\ \text{and that}\ \frac{1}{C_2}<\frac{r_3}{r_2}<C_2.
\end{equation}

%\begin{eqnarray*}
%r_1\leq r_2\ll r_3~~ \mbox{or}~~r_1\ll r_2\approx r_3.
%\end{eqnarray*}
%Here $r_2\ll r_3$ should be interpreted as follows: 
 %For arbitrary large $M>0,$ there exists a convex body $K$ such that its cone volume measure has density $f_K$ %satisfying 
 %\begin{equation}\label{case1}
% 1/\lambda<f_K<\lambda,\ \text{and that}\ \frac{r_3}{r_2}>M.
 %\end{equation}
 %In the following the notation ``$\ll$" is always understood in a similar way. 

We divide our argument into two cases correspondingly.

\begin{figure}[H]
	\centering
	\includegraphics[width=0.45\textwidth]{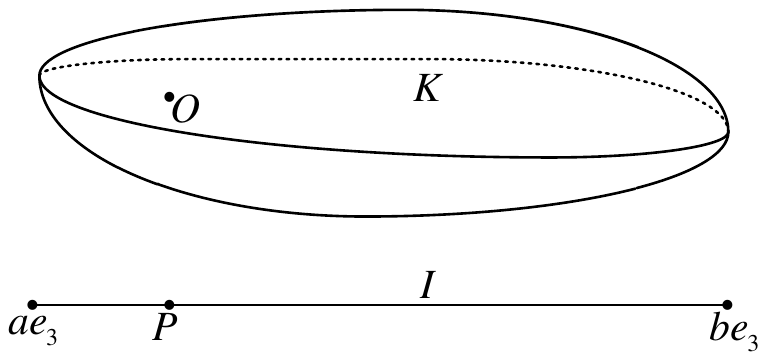}
	\caption{}
    \label{f1}
\end{figure}

{ \noindent\bf CASE I.}
%For any $M$ large, we can find a convex body $K$ satisfying \eqref{case1}.
Denote by $I$ the projection of $K$ on the $e_3$-axis, and $P$ the projection of $O$ (see Fig.~\ref{f1}).

\renewcommand{\figurename}{Fig.}
\renewcommand{\captionlabeldelim}{}

\noindent We now consider two subcases.\\
Subcase i: there exists a constant $\frac{1}{3}>c_0>0,$ such that 
for any $M>0,$ there exists a convex body $K$ satisfying \eqref{case1} and 
\begin{equation}\label{sub111}
\mbox{dist} (P, \partial I)\geq c_0 r_3;
\end{equation}

\noindent Subcase ii: for any $M>0,$ there exists a convex body $K$ satisfying \eqref{case1} and 
\begin{equation}\label{sub22}
\frac{r_3}{\mbox{dist} (P, \partial I)}>M.
\end{equation}

%\begin{equation}\label{sub1}
%\mbox{i:}~~ \mbox{dist} (P, \partial I)\gtrsim  r_3;
%\end{equation}
%\begin{equation}\label{sub2}
%\mbox{ii:}~~ \mbox{dist} (P, \partial I)\ll r_3.
%\end{equation}
%Note that here the subcase i should be understood as follows: there exists a constant $C_0>0,$ such that 
%for any $M>0,$ there exists a convex body $K$ satisfying \eqref{case1} and 
%\begin{equation}\label{sub11}
%\mbox{dist} (P, \partial I)\geq C_0 r_3;
%\end{equation}
%subcase ii should be understood in the following way: for any $M>0,$ there exists a convex body $K$ satisfying \eqref{case1} and 
%\begin{equation}\label{sub22}
%\frac{r_3}{\mbox{dist} (P, \partial I)}>M.
%\end{equation}

\begin{figure}[H]
	\centering
	\includegraphics[width=0.6\textwidth]{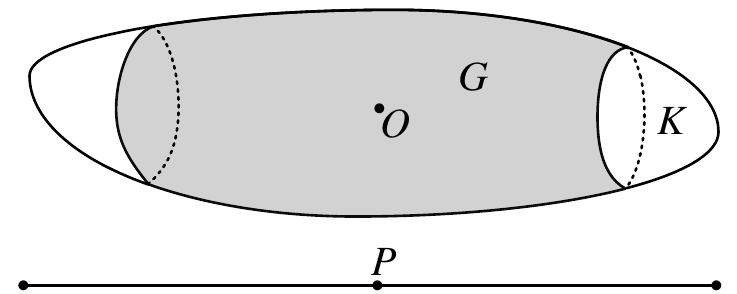}
	\caption{}
    \label{f2}
\end{figure}

{\bf Subcase i.}
Given any $M>0$ large, let $K$ be a convex body satisfying \eqref{case1} and \eqref{sub111}.
 %Let $\nu_x$ be the outer unit normal of $K$ at $x=(x_1, x_2, x_3)$. 
 Take $\epsilon_0=\frac{1}{3}c_0$ and define
\begin{eqnarray*}
F:=\{\nu\in \mathbb{S}^2: \nu\in\nu_K(x)\ \text{for some}\ x\in \partial K, \mbox{dist} (x_3e_3, \partial I)\geq \epsilon_0r_3\},
\end{eqnarray*}
\begin{eqnarray*}
G:=\{x\in \partial K: \mbox{dist} (x_3e_3, \partial I)\geq \epsilon_0r_3\},
\end{eqnarray*}
and see Fig.~\ref{f2}. 
%Note that the difference between $\partial K$ and $\partial' K$ has zero $n-1$-Hausdorff measure.
%
\renewcommand{\figurename}{Fig.}
\renewcommand{\captionlabeldelim}{}

On one hand, by convexity we have that $F$ converges to the equator
\begin{eqnarray*}
\{u\in \mathbb{S}^2: u\cdot e_3=0\},
\end{eqnarray*}
as $\frac{r_2}{r_3}\rightarrow 0$ (choose $M\rightarrow \infty$ in \eqref{case1}). Thus,
\begin{eqnarray*}
\int_Ff\rightarrow 0,
\end{eqnarray*}
as $M\rightarrow \infty.$

On the other hand, by convexity we get
\begin{eqnarray*}
V_K(F)&=&\mbox{volume of the cone with vertex} ~O ~\mbox{and base}~ G\\
&\geq& cr_1r_2r_3\\
&\gtrsim&c,
\end{eqnarray*}
provided $\frac{r_3}{r_2}$ is large enough (guaranteed by choosing $M$ large enough),
where the constant  $c>0$ depends only on $c_0$ and $\lambda.$
This is a contradiction when $M$ is sufficiently large, since
\begin{eqnarray*}
V_K(F)=\frac{1}{3}\int_Ff.
\end{eqnarray*}

{\bf Subcase ii.}
Given any $M>0$ large, let $K$ be a convex body satisfying \eqref{case1} and \eqref{sub22}.
 Let
\begin{eqnarray*}
F:=\{\nu\in \mathbb{S}^2: \nu\in\nu_K(x)\ \text{for some}\ x\in \partial K\backslash B_{\epsilon_0r_3}(o)\}~~\mbox{and}~~G:=\{x:x\in \partial K\backslash B_{\epsilon_0r_3}(o)\}
\end{eqnarray*}
where $\epsilon_0<\frac{1}{2}$ is a small positive constant independent of $M$ and to be determined later.
\renewcommand{\figurename}{Fig.}
\renewcommand{\captionlabeldelim}{}
\begin{figure}[H]
	\centering
	\includegraphics[width=0.6\textwidth]{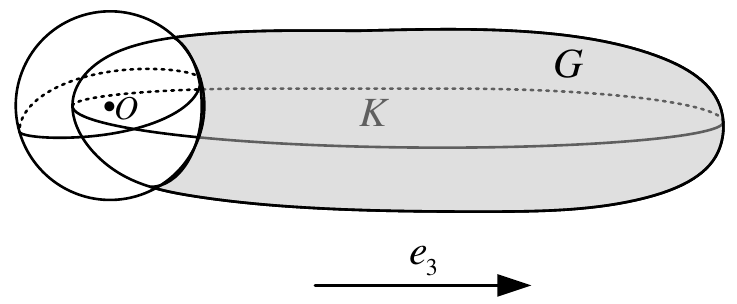}
	\caption{}
    \label{f3}
\end{figure}

On one hand, it follows from convexity that (see Fig.~\ref{f3})
\begin{eqnarray*}
F\rightarrow \mathbb{S}_+^2:=\{u\in \mathbb{S}^2: u\cdot e_3\geq 0\},
\end{eqnarray*}
as $\frac{r_2}{r_3}\rightarrow 0$ (guaranteed by $M\rightarrow \infty$). Hence,
\begin{eqnarray*}
\int_Ff\leq\int_{\mathbb{S}_{\frac{1}{4}}}f=\int_{\mathbb{S}^2}f-\int_{\mathbb{S}^2\backslash \mathbb{S}_{\frac{1}{4}}}f,
\end{eqnarray*}
provided $M$ is sufficiently large,
where
\begin{eqnarray*}
\mathbb{S}_{\frac{1}{4}}:=\{u\in \mathbb{S}^2: u\cdot e_3>-{\frac{1}{4}}\}.
\end{eqnarray*}
Then
\begin{eqnarray*}
\frac{1}{3}\int_Ff\leq \frac{1}{3}\int_{\mathbb{S}^2}f-C_1=V(K)-C_1.
\end{eqnarray*}
for  $C_1=\frac{1}{3}\int_{\mathbb{S}^2\backslash \mathbb{S}_{\frac{1}{4}}}f$ depending only on $\lambda$, provided $M$ is sufficiently large and hence
 $\frac{r_2}{r_3}$ is sufficiently small. 

On the other hand,
\begin{eqnarray*}
V_K(F)&=&\mbox{volume of the cone with vertex} ~O ~\mbox{and base}~ G\\
&\geq&V(K)-\frac{4}{3}\pi\epsilon_0^3.
\end{eqnarray*}
Hence,
\begin{eqnarray*}
V(K)-C_1\geq \frac{1}{3}\int_Ff=V_K(F)\geq V(K)-\frac{4}{3}\pi\epsilon_0^3,
\end{eqnarray*}
which is impossible if initially we choose $\epsilon_0$ small enough such that $\frac{4}{3}\pi\epsilon_0^3<C_1.$
\vskip 10pt

{ \noindent\bf CASE II.} 
In this case, recall that there exists a constant $C_2>0,$ and for any $M>1$ large, there exists a convex body $K$ such that 
\begin{equation}\label{case222}
\frac{r_2}{r_1}>M\ \text{and}\ \frac{1}{C_2}<\frac{r_3}{r_2}<C_2.
\end{equation}
This case is more delicate. Denote by $K^\prime$ the projection of $K$ on the $e_2e_3$ plane, $O^\prime$ the projection of $O$. We first show
\begin{eqnarray}\label{3.1}
\mbox{dist} (O^\prime, \partial K^\prime)\leq \frac{1}{\sqrt{M}C_2}r_3,
\end{eqnarray}
provided $M$ is large enough.
Suppose not, namely for any large $M,$ there exists a convex body $K$ such that \eqref{case222} holds and 
$\mbox{dist} (O^\prime, \partial K^\prime)> \frac{1}{\sqrt{M}C_2}r_3.$
Let
\begin{eqnarray*}
G:=\{x=(x_1, x_2, x_3)\in \partial K:~ \mbox{dist} (x^\prime, \partial K^\prime)\geq \frac{1}{C_3}r_3,~ x^\prime=(x_1, x_2)\}
\end{eqnarray*}
and
\begin{eqnarray*}
F:=\{\nu\in \mathbb{S}^2: \nu\in\nu_K(x)\ \text{for some}\ x\in G\},
\end{eqnarray*}
where $C_3=4\sqrt{M}C_2.$
By convexity, we have
\begin{eqnarray}\label{case24}
\int_{F}f\lesssim |F|\leq C\left(\frac{r_1}{r_3}\right)^2\leq \frac{C}{M^2}
\end{eqnarray}
for some constant $C$ independent of $M,$ provided $M$ is sufficiently large. 
%Hence,
%\begin{eqnarray*}
%\int_{F}f\lesssim r_1^{\frac{4}{3}}.
%\end{eqnarray*}
Moreover,
\begin{eqnarray*}
V_K(F)&=&\mbox{volume of the cone with vertex} ~O ~\mbox{and base}~ G\\
&\geq& \tilde{C}r_1(\frac{r_3}{C_3})^2\geq \frac{\tilde{C}_1}{M} r_1r_2r_3\geq \frac{\tilde{C}_2}{M},
\end{eqnarray*}
where the constants $\tilde{C}, \tilde{C}_1, \tilde{C}_2$ are independent of $M.$
The above estimate contradicts to \eqref{case24} when $M$ is sufficiently large.

\begin{figure}[H]
	\centering
	\includegraphics[width=0.6\textwidth]{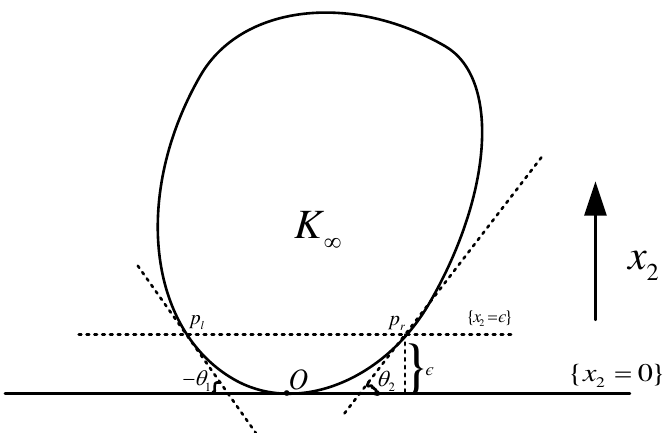}
	\caption{}
    \label{f4}
\end{figure}

Now, we blow down $K$ as follows:
\begin{eqnarray*}
K\longmapsto AK:=\frac{1}{r_3}K.
\end{eqnarray*}
%where
%\[
%A
%=
%\left(
%\begin{array}{ccc}
%\frac{1}{r_3} & 0 &0\\
%0& \frac{1}{r_3}&0\\
%0&0&\frac{1}{r_3}\\
%\end{array}
%\right).
%\]
Then by the assumption on CASE II: $r_1\ll r_2\approx r_3$, namely \eqref{case222} holds,
it follows that
\begin{eqnarray*}
AK=\frac{1}{r_3}K\rightarrow K_\infty\subset \mathbb{R}^2
\end{eqnarray*}
in Hausdorff distance, as $M\rightarrow \infty.$
Moreover,
\begin{eqnarray*}
B_{\frac{1}{C_3}}(z)\subset K_\infty\subset B_{C_3}(z)
\end{eqnarray*}
for $z\in K_\infty$ and some constant $C_3>0$. By (\ref{3.1}), we see $o\in \partial K_\infty$. By a rotation of coordinates we may assume
\begin{eqnarray*}
K_\infty\subset \{x_2\geq 0\},
\end{eqnarray*}
and see Fig.~\ref{f4}.
\renewcommand{\figurename}{Fig.}
\renewcommand{\captionlabeldelim}{}
\noindent Keep in mind that $\frac{1}{r_3}K$ can be as close to $K_\infty$ as we want, provided $M$ is large enough. 
%Then there exist two subcases:

%{\noindent\bf Subcase i:} $K_\infty$ is strictly convex at $o$, namely
%\begin{eqnarray*}
%\{x_2=0\} \cap \partial K_\infty=\{o\};
%\end{eqnarray*}

%{\noindent\bf Subcase ii:} $\{x_2=0\} \cap \partial K_\infty$ contains a segment.

Recall the spherical coordinate for $\mathbb{S}^2$:
\begin{eqnarray*}
T:~(\theta, \varphi)\in(-\pi, \pi)\times \left(-\frac{\pi}{2}, \frac{\pi}{2}\right)\longrightarrow (\cos \varphi\cos\theta, \cos\varphi\sin\theta, \sin\varphi)\in \mathbb{S}^2.
\end{eqnarray*}
Let
\begin{eqnarray}\label{3.2}
D:=\{x_2\leq \epsilon\}\cap \partial K_\infty
\end{eqnarray}
and
\begin{eqnarray*}
\widetilde{F}:=\{\nu\in \mathbb{S}^1: \nu\in\nu_{K_{_\infty}}(x)\ \text{for some}\ x\in D\},
\end{eqnarray*}
%where $\widetilde{\nu}_x$ denotes the outer unit normal of $K_\infty$ (as a two dimensional object) at $x$. 
Let $p_l, p_r$ be the intersections of the line $\{x_2=\epsilon\}$ and $\partial K_\infty.$ Since $\partial' K_\infty$ is dense in $\partial K_\infty,$ we can find a sequence of point $z_k\in \partial' K_\infty\cap D$ converges to $p_r,$ where $\partial' K_\infty$ is the subset of 
$\partial K_\infty$ where $\partial K_\infty$ is differentiable. Note that $\partial K_\infty\backslash \partial 'K_\infty$ has $1$-Hausdorff measure 0.
 Let $\alpha_k$ be the cute angle between the tangent line of $\partial K_\infty$ at $z_k$ and the $x_1$-axis.
 Then by convexity $\alpha_k\rightarrow \theta_2$ as $k\rightarrow \infty$ for some $\theta_2\in (0, \pi/2),$ see Fig.~\ref{f4}. $\theta_1$ is defined similarly.

By convexity, we have
\begin{eqnarray*}
(-\theta_1, \theta_2)\subset \widetilde{F}\subset \mathbb{S}^1.
\end{eqnarray*}
Observe that $\frac{\theta_1+\theta_2}{\epsilon}\rightarrow \infty$ as $\epsilon\rightarrow 0.$
%where $\theta_1+\theta_2\gg \epsilon$.
%We first rule out subcase i.
 Let
\begin{eqnarray*}
F:=T\left(\left(-\frac{\theta_1}{2}, \frac{\theta_2}{2}\right)\times \left(-\frac{\pi}{4}, \frac{\pi}{4}\right)\right).
\end{eqnarray*}
For any $\xi=T(\theta, \alpha)\in F$, there exists a point $z\in D$ such that $T(\theta, 0)\in \nu_{K_{_\infty}}(z)$. Thus
\begin{eqnarray}\label{Hdef}
H:=\{x: ~\xi\cdot (x-z)=0\}
\end{eqnarray}
passes $z$, and
\begin{eqnarray} \label{move1}
\mbox{dist}(H, K_\infty\cap \{x_2>2\epsilon\})>\delta>0
\end{eqnarray}
for some positive $\delta$ depending only on $\partial K_\infty$, $\epsilon$, $\theta_{_1}$ and $\theta_{_2}$. Since $\frac{1}{r_3}K\rightarrow K_\infty$ in the Hausdorff distance, as $M\rightarrow \infty$, we have that if $\nu_K(y)\cap F\ne \emptyset$ for some $y\in \partial \left(\frac{1}{r_3}K\right),$ then
\begin{eqnarray}\label{3.3}
y\in\{x_2\leq 2\epsilon\}\cap \partial \left(\frac{1}{r_3}K\right).
\end{eqnarray}
Indeed, letting $H$ be as in \eqref{Hdef} and $\xi\in \nu_K(y)\cap F$, it follows  from \eqref{move1} and the fact that 
$\frac{1}{r_3} K$ converges to $K_\infty$ that 
\begin{eqnarray} \label{move2}
\mbox{dist}(H, \frac{1}{r_3}K\cap \{x_2>2\epsilon\})>\frac{1}{2}\delta
\end{eqnarray}
provided $M$ is sufficiently large.
Now, we translate $H$ such that it becomes a supporting plane of $\partial \left(\frac{1}{r_3} K\right)$, then the touching point must be inside $\{x_2\leq 2\epsilon\}$. Hence \eqref{3.3} holds, which implies that
 $\nu_K^{-1}(F)\subset \{x_2\leq 2\epsilon\}.$

 By the definition of $F$ we have
\begin{eqnarray*}
\int_Ff\geq C(\theta_1+\theta_2),
\end{eqnarray*}
for some $C$ depending only on $\lambda,$ provided $M$ is sufficiently large.
Moreover,
\begin{eqnarray*}
V_K(F)&=&\mbox{volume of the cone with vertex} ~O ~\mbox{and base}~ \{x_2\leq 2\epsilon r_3\}\cap \partial K\\
&\leq& C_1\epsilon r_1r_2r_3\\
&\leq& C_2\epsilon,
\end{eqnarray*}
where $C_1, C_2$ depending only on $\lambda,$ provided $M$ is sufficiently large.
This is a contradiction for $M$ sufficiently large, since
\begin{eqnarray*}
V_K(F)=\frac{1}{3}\int_Ff~~\mbox{and}~~ \frac{\theta_1+\theta_2}{\epsilon}\rightarrow \infty,
\end{eqnarray*}
as $M\rightarrow \infty.$ 

%Now we give the proof of subcase ii. Let $\eta=\left\{\nu_y:  y\in\{x_2=0\} \cap \partial K_\infty\right\}$. We take
%\begin{eqnarray*}
%F=T\left(\eta\times \left(-\frac{\pi}{4}, \frac{\pi}{4}\right)\right).
%\end{eqnarray*}
%Hence, for any $\xi\in F$, there exists a point $z\in D$ such that $\widetilde{\nu}_z=T(\theta, 0)$. 
%Obviously,
%\begin{eqnarray*}
%H:=\{x: ~\xi\cdot (x-z)=0\}
%\end{eqnarray*}
%passes the point $z$, and
%\begin{eqnarray*}
%\mbox{dist}(H, K_\infty\cap \{x_2>\epsilon\})>\delta>0
%\end{eqnarray*}
%for some positive $\delta$ depending only on $\partial K_\infty$, $\epsilon$, $\theta_{_1}$ and $\theta_{_2}$. Note that $\frac{1}{r_3}K\rightarrow K_\infty$, as $\frac{r_1}{r_3}\rightarrow 0$. Similar to the derivation of (\ref{3.3}), if $y\in \partial \left((1/r_3)K\right)$ and $\nu_y\in F$ then
%\begin{eqnarray*}
%y\in\{x_2\leq \epsilon\}\cap \partial \left(\frac{1}{r_3}K\right).
%\end{eqnarray*}
%Hence,
%\begin{eqnarray*}
%\theta_1+\theta_2\lesssim V_K(F)&=&\mbox{volume of the cone with vertex} ~O ~\mbox{and base}~ \{x_2\leq %\epsilon r_3\}\cap \partial K\\
%&\lesssim&\epsilon r_1r_2r_3\\
%&\lesssim&\epsilon,
%\end{eqnarray*}
%which contradicts
%\begin{eqnarray*}
%\theta_1+\theta_2\gg \epsilon.
%\end{eqnarray*}

\section{\bf Proof of Theorem 1.1}

\ \ \ \ \  Let $K\in \mathcal{K}^n$, satisfying
\begin{eqnarray}\label{4.1}
\mbox{det}(\nabla^2h_K+h_KI)=\frac{f}{h_K}.
\end{eqnarray}
We compute the linearized equation of (\ref{4.1}) at $h_K=1$ as follows:
\begin{eqnarray*}
U^{ij}(\varphi_{ij }+\varphi\delta_{ij})=-h_K^{-2}f\varphi,
\end{eqnarray*}
where $U^{ij}$ is the cofactor matrix of $\nabla^2h_K+h_KI$. Since $h_K=1$, we have
\begin{eqnarray*}
\triangle_{\mathbb{S}^2}\varphi+3\varphi=0,
\end{eqnarray*}
which is invertible. Thus by the inverse function theorem and Schauder estimate we have

{\it \noindent{\bf Lemma 4.1.}~There exists a small constant $\epsilon_0>0$ such that if $\|f-1\|_{C^\alpha}<\epsilon_0$, $\|h_K-1\|_{L^\infty}\leq \epsilon_0$ and $\|h_L-1\|_{L^\infty}\leq \epsilon_0$ where $K, L\in \mathcal{K}^3$ satisfy equation (\ref{1.2}), then $K=L$.

}
{\bf Proof.} By Caffarelli's $C^{2, \alpha}$ estimate for Monge-Amp\`ere equation (see \cite{C901, C902}), we have
\begin{eqnarray*}
\|h_K\|_{C^{2, \alpha}}\leq C_0~~\mbox{and}~~\|h_L\|_{C^{2, \alpha}}\leq C_0,
\end{eqnarray*}
for some constant $C_0>0$ depending only on $\epsilon_0.$ Thus,
\begin{eqnarray*}
\left\|\frac{f}{h_K}-1\right\|_{C^{\alpha}}\leq \widetilde{\epsilon}_0\rightarrow 0 ~~\mbox{as}~~\epsilon_0\rightarrow 0,
\end{eqnarray*}
and
\begin{eqnarray*}
\left\|\frac{f}{h_L}-1\right\|_{C^{\alpha}}\leq \widetilde{\epsilon}_0\rightarrow 0 ~~\mbox{as}~~\epsilon_0\rightarrow 0.
\end{eqnarray*}
We compute
\begin{eqnarray}\label{4.2}
\frac{f}{h_K}-1&=&\mbox{det}(\nabla^2h_K+h_K\delta_{ij})-\mbox{det}(\nabla^2 1+\delta_{ij}) \nonumber \\
&=&\int_0^1\frac{d}{dt}\mbox{det}\left(\nabla^2((1-t)+th_K)+((1-t)+th_K)I\right)dt \nonumber \\
&=&\sum_{i,j=1}^{n-1}\int_0^1U_t^{ij}dt\cdot\left((h_K-1)_{ij}+(h_K-1)\delta_{ij}\right)\nonumber \\
&=&\sum_{i,j=1}^{n-1}a_{ij}\left((h_K-1)_{ij}+(h_K-1)\delta_{ij}\right),
\end{eqnarray}
where the coefficient $a_{ij}=\int_0^1U_t^{ij}dt$, and $U_t^{ij}$ is the cofactor matrix of
\begin{eqnarray*}
\nabla^2((1-t)+th_K)+((1-t)+th_K)I.
\end{eqnarray*}
Since
\begin{eqnarray*}
\|h_K\|_{C^{2, \alpha}}\leq C_0,
\end{eqnarray*}
we have
\begin{eqnarray*}
\frac{1}{C_1}I\leq \{a_{ij}\}\leq C_1I,
\end{eqnarray*}
for some constant $C_1>0$. Namely, (\ref{4.2}) is uniformly
elliptic. Hence, by Schauder estimate we have
\begin{eqnarray*}
\|h_K-1\|_{C^{2, \alpha}}&\leq &C\left(\|h_K-1\|_{L^\infty}+\left\|\frac{f}{h_K}-1\right\|_{C^{\alpha}}\right)\\
&\leq &C(\epsilon_0+\widetilde{\epsilon}_0).
\end{eqnarray*}
Similarly,
\begin{eqnarray*}
\|h_L-1\|_{C^{2, \alpha}}\leq C(\epsilon_0+\widetilde{\epsilon}_0).
\end{eqnarray*}
Therefore, by the inverse function theorem we have
\begin{eqnarray*}
K=L,
\end{eqnarray*}
provided $\epsilon_0$ is sufficiently small.

Now we can use compact argument to finish the proof of Theorem 1.1.

{\bf Proof of Theorem 1.1.}~~If not, then by Lemma 4.1 we have that there exist $f_i$, $K_i\in \mathcal{K}^n$ such that
\begin{eqnarray*}
{h_{K_i}}\mbox{det}(\nabla^2h_{K_i}+h_{K_i}I)=f_i,
\end{eqnarray*}
and when $i\rightarrow \infty$,
\begin{eqnarray*}
\|h_{K_i}-1\|_{L^\infty}>\epsilon_0,~~ \|f_{i}-1\|_{C^\alpha}\rightarrow 0.
\end{eqnarray*}
By Lemma 1, we see
\begin{eqnarray*}
\|h_{K_i}\|_{L^\infty}\leq C
\end{eqnarray*}
for some constant $C$ depending only on $\lambda.$
Hence, by Blaschke's selection theorem, up to a subsequence we may assume $K_i\rightarrow K_\infty$ in Hausdorff distance. Thus by weak convergence of surface area measure we have
\begin{eqnarray*}
{h_{K_\infty}}\mbox{det}(\nabla^2h_{K_\infty}+h_{K_\infty}I)=1,
\end{eqnarray*}
where
\begin{eqnarray*}
\|h_{K_\infty}-1\|_{L^\infty}>\epsilon_0.
\end{eqnarray*}
Hence \eqref{1.2} admits two different solutions $h_{K_\infty}$ and $1,$  
which contradicts to the uniqueness of solutions to the logarithmic Minkowski problem with constant $f$ (see \cite{48a, BCD}).

\end{document}